\documentclass[12pt,a4paper]{article}

\usepackage{amsmath}
\usepackage{amssymb}
\usepackage{amsxtra}
\usepackage{latexsym}
\usepackage{color}

\newtheorem{pr}{Proposition}

\newtheorem{lem}{Lemma}

\newtheorem{ex}{Example}
\newtheorem{as}{Assumption}
\newtheorem{rem}{Remark}
\newtheorem{de}{Definition}

\let\<\langle
\let\>\rangle

\numberwithin{equation}{section}  

\begin{document}

\title{On differential operators and linear differential equations on torus}

\author{
Vladimir P. Burskii
\thanks{Moscow Institute of Physics and Technology,
{\tt bvp30@mail.ru}}
}

\date{}

\maketitle

\begin{abstract}

In this paper, we consider periodic boundary value problems for differential equations whose coefficients are trigonometric polynomials.
We construct the spaces of generalized functions, where such problems have solutions.
In particular, the solvability space of a periodic analogue of the Mizohata equation is constructed.
We build also a periodic analogue and a generalization of the construction of the nonstandard analysis,
where infinitely smalls are not only functions, but also functional spaces.

To show that not all constructions on the torus lead to a simplification in compare with the plane,
we consider a periodic analogue of the hypoelliptic differential operator and show that its number-theoretic properties are significant.
In particular, it turns out that if a polynomial with integer coefficients is irreducible in the rational field,
then the corresponding differential operator is hypoelliptic on the torus.


\end{abstract}

\section{Introduction}

Periodic boundary value problems for differential equations is a famous object both in mathematical education and research (see, for example, [1], [2], [3]).
As Lax noted in [4], in the periodic theory of differential equations is free of some technical difficulties that arise in non-periodic theory.
This makes it possible to create a more beautiful theory.

A study of periodic boundary value problems for linear differential equations brings us to the wonderful world of functions on the torus.
Since the torus as the product of a finite number of circles, dealing with the torus simplifies studying the behavior of functions of several variables by each variable separately. Moreover, the basis elements are eigenfunctions of linear differential operators with constant coefficients.
In addition, the topology of the torus allows us to forget about the boundary and the behavior at infinity.
This allows us to focus our attention on the only infinity that we face up in this way:
the infinite dimension of functional spaces.

There are many works devoted to periodic boundary value problems for differential equations considered as equations on the torus.
See, for example, the books [1], [2], [3] and the references therein.
However, the present paper has no intersection with them.
A viewpoint to the structure of the set of infinitesimals presented in this paper is essentially different from those in the non-standard analysis
(see, for example, [14]). The first mention of the presented results was published in Russian in [15],[16]. In the future it would be interesting to consider, in particular, the spectral properties of differential operators on a torus in the spirit of elliptic theory.

\section{Spaces of periodic functions}

\subsection{Spaces of periodic functions}

The number of variables in this problem is not significant,
and without loss of generality we shall consider the case of two variables.

It is well known that every Fourier series (trigonometric series)
$$
\sum\limits_{k \in {\Bbb Z} \oplus {\Bbb Z}} a_k e^{i k x}, \ \ \sum_k |a_k|^2 < \infty,
$$
where $a_k \in {\Bbb C}$, $k = (k_1, k_2)$, $x = (x_1, x_2) \in {\Bbb R}^2$, and $k x = k_1 x_1 + k_1 x_2$,
defines a periodic square integrable function on $x$.
Thus, such series form the space $L_2(T^2)$, where $T^2={\Bbb R}^2/{\Bbb Z}^2$ is the torus of dimension~2.

For every $m \in {\Bbb Z}$, $m \ge 2$, define by $H^m$ the space of $2\pi$-periodic complex-valued functions in ${\Bbb R}^2$ such that
$$
\Vert u \Vert_m^2 = \int\limits_{T^2} \overline u (x) (1 - \Delta)^m u (x) d x < \infty,
$$
where $\Delta = \frac {\partial^2}{\partial x_1^2} + \frac {\partial^2}{\partial x_2^2}$ is the Laplace operator.
All $H^m$ are Hilbert spaces, namely, the famous Sobolev spaces.

It is known (see, for example, [1]) that functions $\exp ( i k x)$, where $k x = k_1 x_1 + k_2 x_2$,
form an orthogonal basis in $H^m$ and, consequently, every function $f \in H^m$ is expandable into the Fourier series
$$
f = \sum_{k \in {\Bbb Z\oplus Z}} f_k e^{i k x}
$$
converging to $f$ in the topology of $H^m$.

Further we shall consider trigonometric series with arbitrary real coefficients, not necessarily converging (such series are usually called formal).
One can consider $H^m$ as the vector space consisting of formal Fourier series with the finite norm
$$
\Vert f \Vert^2_m = \sum_k (1 + k \cdot k)^m \cdot |f_k|^2.
$$
This formula defines a norm in the space $H^m$ with any real $m$.

Consider the vector space $H^m$
with the topology of the space ${\Bbb R}^{\Bbb Z}$, in which $H^m$ is continuously embedded.
For $m < 0$, the space $H^m$ is conjugate to the space $H^{-m}$
in the topology of the space $H^0 = L_2 (T^2)$.


Moreover, can prove a more general statement:

\begin{pr}
\label{St1}
Let $E$ be a barreled vector topological space of $2\pi$-periodic functions continuously embedded in $H^0$.
Let the system $\{ e^{i k x} \}$ be a basis in $E$, that is,
for every $v \in E$ there exists a unique sequence $\{ v_k \} \subset E$ such that
$$
\sum_{k^2 \le N} v_k e^{i k x} \to v  \ \ \text{as} \ \ N \to \infty.
$$
Then the dual space $E^*$ is naturally isomorphic to a subspace of the space
$$
F = \bigl\{  u = \sum_k u_k e^{i k x} : \
\sum v_k e^{i k x} \in E, \  \< u, v\> = \sum u_k \overline v_k < \infty \ \
\forall v \in E
\bigr\}.
$$
Here the paring $\< \cdot, \cdot \>$ gives rise to the duality.
\end{pr}


{\bf Proof}.
This statement follows from the known fact that the Mackey topology in barrel spaces coincides with the original topology,
since it is the strongest among all topologies consistent with duality (see [5]).

\medskip

There are several examples, which we shall use below:

\begin{ex}
\label{Ex1}
{\rm
The space of infinitely differentiable periodic functions $H^\infty = \underset {m}{\cap} H^m$, whose total element has the form
$$
u = \sum u_k e^{i k x}, \ \  k^{2l} u_k \overset {\scriptsize k \to \infty}{\longrightarrow} 0 \ \ (\forall \, l).
$$
Also the conjugate space $(H^\infty)^* = H^{-\infty}= \underset {m}{\cup} H^m$, which is the space of periodic distributions whose Fourier coefficients
tend to infinity no faster than some power of $k^{2l}$.
}
\end{ex}

\begin{ex}
\label{Ex2}
{\rm
The space $E_0$ of $2\pi$-periodic functions
$$
u = \sum u_k e^{i k x}, \ \
\exists \delta_1 >0, \ \delta_2 > 0: \ \, \sum e^{|k_1| \delta_1 + |k_2| \delta_2} |u_k| < \infty,
$$
included in the space of periodic real analytic functions.
Another examples is the conjugate space $E_0^*$, which consists of series $\sum u_k e^{i k x}$ whose coefficients $u_k \overset {\scriptsize k \to \infty}{\longrightarrow} \infty$ slower than any exponential $e^{|k_1| \delta_1 + |k_2| \delta_2}$. This space contains the space of hyperfunctions ([6]).
}
\end{ex}

\begin{ex}
\label{Ex3}
{\rm
The space $l_1 (|k_1 |!)$ of functions
$$
u = \sum u_k e^{i k x}, \ \ \sum_k |k_1 |! \, |u_k | < \infty.
$$
Also the conjugate space  $l_1^* (|k_1| !)$ of series
$\sum_k v_k e^{i k x}$, $v_k = O (|k_1|!)$.
}
\end{ex}

Let $P (x_1, x_2)$ be a homogeneous polynomial of degree $p$ with constant coefficients.
Consider the differential operator $\hat P : H^m \to H^{m - p}$ generated by the polynomial $P$:
$$
\hat P u = P \biggl( - i \frac{\partial}{\partial x_1}, - i \frac {\partial}{\partial x_2} \biggr) u.
$$
The obvious formula
$$
\hat P \biggl( \sum_{k \in {\Bbb Z\oplus Z}} f_k e^{i k x} \biggr) = \sum_{k \in {\Bbb Z\oplus Z}} P(k_1,k_2)f_k e^{i k x}
$$
allows us to consider the operator $\hat P$ on the space $F$ of formal trigonometric series.

\subsection{Solvability of the Mizohata equation}

In [7],  G.~Levy gave an example of a linear differential equation of the first order with infinitely differentiable coefficients
that has no solutions in the space of distributions in three-dimensional space.
Developing the ideas of G.~Levy and P.~Garabedyan [9], V.\,V.~Grushin in [8] gave an example of a first order differential equation with infinity differentiable coefficients
that has no solutions in the space of distributions on the plane:
\begin{equation}
\frac {\partial u}{\partial x} + i x \frac {\partial u}{\partial y} = f(x, y).
\label{(1)}
\end{equation}
The operator in the left-hand side of equation \eqref{(1)} is one of the Mizohata operators, considered in [10].
The function $f \in C_0^\infty ({\Bbb R}^2)$ is even by~$x$, it was constructed by Grushin in a special way.

Consider a periodic modification of equation \eqref{(1)}:
\begin{equation}
\frac{\partial u}{\partial x} + i \sin{x} \, \frac{\partial u}{\partial y} = \tilde f (x,y),
\label{(2)}
\end{equation}
where $\tilde f$ is $2\pi$-periodic continuation of the function $f$ mentioned above.
It can be checked that the Grushin's reasonings are also applicable to equation~\eqref{(2)}.

We shall prove that equation \eqref{(2)} has a solution in a wider space of generalized functions than the space of Schwartz distributions.

\begin{pr}
\label{St2}
For any even right-hand side $\tilde f \in H^{-\infty}$ equation \eqref{(2)} has a unique periodic solution $u(x_1, x_2)$ odd in the variable $x_1$,
which belongs to the space $l_1^* (|k_1| !)$.
\end{pr}

{\bf Proof}.
Let us write the equation \eqref{(2)} in the form
$$
\frac {\partial u}{\partial x_1} + \frac {e^{i x_1} - e^{- i x_1}}2 \frac {\partial u}{\partial x_2} \ =\tilde f,
$$
which yields
\begin{equation}
k_1 u_{k_1, k_2} + \frac {k_2}2 \left ( u_{k_1 -1, k_2} - u_{k_1 + 1, k_2} \right ) = f_k.
\label{(3)}
\end{equation}
For a fixed $k_2 \not= 0$ we obtain the recurrent formula with respect to $k_1$
\begin{equation}
u_{k_1 + 1, k_2} = \frac 2{k_2} \left ( k_1 u_{k_1, k_2} - f_k \right) u_{k_2 - 1, k_2}.
\label{(35)}
\end{equation}

Since the function $u(x_1, x_2)$ is odd in the variable $x_1$, we have $u_{0, k_2} = 0$, $u_{-k_1, k_2} = - u_{k_1,k_2}$.
Therefore, the coefficients $u_k$ are uniquely determined by \eqref{(35)}.
From \eqref{(35)} we have the following estimation:
$$
|u_{k_1 + 1, k_2}| < \sum_{j  = 0} (j + 1) ! f_{k_1 - j} < (k_1 + 1) ! \sum_k f_k = c (k_1 + 1) !
$$
Thus, the solution $u = \sum u_k e^{i k x}$ belongs to the space $l_1^* (|k_1| !)$.

\medskip

Note that the operations of differentiation and multiplication by a trigonometric polynomial defined formally in the space $F$,
coincide with the analogues operations in the Banach space $l_1^* (|k_1| !)$ defined as usual in spaces of generalized functions through pairing.

\begin{pr}
\label{St3}
Every periodic solution $u (x_1, x_2)$ of homogeneous equation \eqref{(2)} is even in $x_1$ and it is uniquely determined by the functions
$$
u_0 (x_2) : = \int_0^{2 \pi} u(x_1, x_2) d x_1 : = \<u, 1 \>_{x_1},
$$
$$
u_1 (x_2) : = \int_0^{2 \pi} u (x_1, x_2) e^{- i x_1} d x_1 : = \< u, e^{i x_1} \>_{x_1}.
$$
It belongs to the space $l_1^* (|k_1| !)$ if $u_0$ and $u_1$ have bounded sequences of coefficients.
\end{pr}

{\bf The proof} follows from formula~\eqref{(3)}.

Here we consider the function $u$ as a formal trigonometric series, and pairing along one coordinate is defined in the standard way:
$$
\<u, v\>_{x_1} : = \sum_n \Bigl\< \sum_k u_{kn} e^{i k x_1}, \ \sum_m v_{mn} e^{i m x_1} \Bigr\> e^{inx_2}=
$$
$$
= \sum_n \Bigl( \sum_m u_{m\,n} \overline v_{-m\,n} \Bigr) e^{i n x_2}.
$$
It is clear that pairing on $x_1$ does not always exist, but if $t$ is a trigonometric polynomial, then the function $\<t, v\>_{x_1} (x_2)$ exists.

\subsection{Solvability of general equations}

Now let us consider the general operator
$L : \sum_{|\alpha| \le m} T_\alpha (x) D^\alpha$, where
$T_\alpha (x)$ is a trigonometric polynomial of degree $(s_\alpha^1, s_\alpha^2)$.
It can also be written in the form
$$
L = \sum^{s^1}_{n_1 = - s^1} \sum^{s^2}_{n_2 = - s^2} e^{i n x} P_n (D), \ \ \
s^1 : = \underset {\alpha}{\max}\ s^1_\alpha, \ \,
s^2 : = \underset {\alpha}{\max}\ s^2_\alpha.
$$

Let us assume that the operator $L$ satisfies the following condition:

\begin{as}
\label{Asum}
For every $n$ the equation $P_n (x) = 0$ has no solutions in integers.
\end{as}


Then the following generalization of Proposition~\ref{St3} is true.


\begin{pr}
\label{St4}
Under Assumption~\ref{Asum}, every formal periodic solution $u (x_1, x_2)$ of equation $L u = 0$
is uniquely determined by the functions
\begin{align*}
u_{02} (x_2) : = \<u, 1\>_{x_1}, \ \
u_{12} (x_2) : = \<u, e^{i x_1} \>_{x_1}, \, \ldots, \,
u_{s^1 2} (x_2) : = \<u, e^{i s^1 x_1} \>_{x_1}, \\
u_{01} (x_1) : = \<u, 1\>_{x_2}, \ \
u_{11} (x_1) : = \<u, e^{i x_2} \>_{x_2},\, \ldots, \,
u_{s^2 1} (x_1) : = \<u, e^{i s^2 x_2} \>_{x_2}.
\end{align*}
if they satisfy the following conditions:
$$
\< u_{q 1}, e^{i p x_1} \> _{x_1} = \< u_{p 2}, e^{i q x_2} \>_{x_2}, \ \ \forall \, p,q = 1, \ldots, \min (s^1, s^2).
$$
\end{pr}

{\bf The proof} is by the direct substitution of formal series into the equation.

\subsection{Linear sections as objects of non-standard analysis}

Consider the vector space of formal trigonometric series $F$ and define a relation of order $u \le v$ in the following way.
An element $u \in F$ is {\it more regular} than $v \in F$ or, equivalently, $v$ is {\it more singular} than $u$ if
there exists an element $h \in F$ with bounded positive coefficients $h_k$ such that $u$ is the convolution of $v$ and $h$:
$$
u = v * h := \sum_k v_k h_k e^{i k x}.
$$
This is obviously equivalent to the condition
$$
\exists C > 0: \ \ \Bigl| \frac {u_n}{v_n} \Bigr| < C \ \  \forall n \in {\Bbb Z}.
$$
Similar order relations are used in asymptotic expansions [11].

We shall the following definition:
a subspace $\alpha \subset F$ is called a linear section of $F$ if from $v \in \alpha$ it follows that $u \in \alpha$ for every $u \le v$.
A trivial example: if $v \in F$, then the set of $u \in F$ such that $u \le v$ is a linear section of $F$.
It is called the principal linear section of $F$.

All the subspaces considered above are also linear sections.
Note that the pairing considered in Section~1 naturally generates a Hausdorff topology on every vector subspace of $F$ (see [5]).
Therefore, every linear section $\alpha$ can be considered as a complete topological vector space,
and the space $\alpha^*$ consisting of $g \in F$ such that $\<f,g\> < \infty$ for all $f \in \alpha$ is the dual space of $\alpha$.

The set $M$ of all linear sections of $F$ is ordered by the inclusion.
For every two sections $\alpha$ and $\beta$ there exist the linear sections
$$
\mbox{sup} (\alpha, \beta) = \alpha +\beta, \ \ \mbox {inf} (\alpha, \beta) = \alpha \cap \beta
$$
defined as minimal (by the inclusion) linear sections that contain respectively $\alpha \cup \beta$ or $\alpha \cap \beta$.

It is easy to see that the distributivity relations are satisfied:
$$
\alpha \cap (\beta + \gamma) = \alpha \cap \beta + \alpha \cap \gamma, \ \
\alpha + (\beta \cap \gamma ) = (\alpha + \beta) \cap (\alpha + \gamma).
$$
Thus, in the set $M$ a certain structure of the distributive lattice with additive and multiplicative identity elements is introduced.
If $A$ is an arbitrary set and $\{\delta_a|\ a \in A\}$ is a family of linear sections,
then the supremum of $\mbox {sup}\ \delta_a$ is minimal linear section that contains all $\delta_a$.
By the Zorn lemma, the supremum of any family exists.

\subsection{Linear sections and solvability of general equations}

\begin{pr}
\label{St5}

1. The operator
$$
L = \sum_{|\alpha| \le m}T_\alpha (x) D^\alpha
$$
sends any linear section to a linear section an, consequently, it induces a mapping (endomorphism) $\tilde L$ in the set $M$, which
preserves the lattice structure.

2. Under Assumption~\ref{Asum},
the mapping $\tilde L$ is an epimorphism of the lattice $M$ and for every linear section $G$ there exists
the maximum $\beta$ among those linear sections $\alpha$ for which $\tilde L \alpha \le G$, and $\beta \neq F$.
\end{pr}

{\bf Proof}.
Since operations of differentiation, multiplication by a scalar, addition, and shifts $\{u_n \}\to \{ u_{n + k}\}$ preserve the relation of order $\le$,
the mapping $L$ sends every linear section to a linear section and the induced mapping preserves the relation of order and the lattice operations.

From Assumption \ref{Asum}, it follows that every equation $L u = e^{i n x}$ has a solution $u_n$ among trigonometric polynomials.
Further, if $f$ is a formal series and the operator $L$ satisfies condition~\ref{Asum}, then the equation
$$
L u = f = \sum_n f_n e^{i n x}
$$
is solvable in the class of formal series. Let such a solution be the formal series
$$
w = \sum_{n, k} f_n u_{n k} e^{i k x},
$$
where for each $k$ the sum over $n$ is finite, and the coefficients $u_{nk}$ are uniquely defined by condition~\ref{Asum}.
Given linear section $G$, for every $f \in G$ consider the set $W_f$ of solutions to the equation $L u = f$ and the set of the corresponding principal linear sections.
By the Zorn lemma, there exists the supremum $s\in \tilde W_f$, which is the desired linear section $\beta$.
Obviously, $\beta$ does not coincide with $F$, since otherwise $G = F$.
This completes the proof.


\begin{rem}
{\rm
Assumption~\ref{Asum} can be replaced with be replaced with a weaker condition:
$$
\forall m \in {\Bbb Z}^2 \ \  \exists n: \ P_n(m) \neq 0.
$$
}
\end{rem}


\begin{de}
The linear section $\beta$ constructed in Proposition~\ref{St5} is called
a solution of the equation $L u = G$ with linear section $G$ in the right-hand side.
\end{de}

For example, from what we proved above, it follows that solution of equation \eqref{(2)} with the right-hand side $G = H^m$ is the section $\beta \subset l_1^* (|n_1| !)$.

\begin{rem}
{\rm

The term ``section'' was chosen due to the obvious analogy with Dedekind sections when constructing the field of real numbers.
Note also that the presented construction is consonant with some constructions of non-standard analysis related to the extension of the field of reals:
the infinitesimal germs of functions from both a field and an ultrafilter.

The set $M$ of linear sections is only partially ordered, but it contains all the germs of sequences as principal sections.
In the set $M$, there exists the operation of convolution, which is associative, commutative and distributive with respect to addition and intersection
(however, the inverse exists not for all elements). Moreover, in $M$ there exists the conjugation and all linear sections are reflexive spaces.
In the set $M$ one can also introduce an associative commutative product, which is not always defined,
but covering the product of smooth functions and the product of distributions according to Mikusinsky--Hirata--Ogawa, and, in addition, the inverse element not lying in $M$ as, so to speak, a singular section, containing not all more regular sequences, but all more singular, and so on. The most important thing here is that this set contains not only
all functions, but also all functional spaces, which, along with each of their elements, also contain increasingly smooth elements. And now we can consider the question of solving the differential equation in a class of function spaces, if the given right-hand side is a space as it was in statement 5.
}
\end{rem}


\section{On the hypoellipticity of differential operators on the torus}

As we noted above, according to Lax's statement ([4]), in the periodic case it is possible to construct a more beautiful theory.
However, the periodic case is not always simpler that the non-periodic one.
In this section, we characterize homogeneous differential operators with constant coefficients hypoelliptic in the space of periodic functions on the plane.
This is one of the cases when the Lax's statement is quite controversial.

Let $H^m$, $m \ge 0$, be the space of complex functions on the plane, $2\pi$-periodic in the both arguments such that
$$
\Vert u \Vert_m^2 = \int\limits_{T^2} \overline u (x) (1 - \Delta)^m u (x) dx < \infty, \ \ \,
\Delta = \frac {\partial^2}{\partial x_1^2} + \frac {\partial^2}{\partial x_2^2}.
$$

The space $H^{- m}$ is defined as the space dual to $H^m$ in the $H^0$-topology.
In Section~1, we noted that the functions
$$
\exp ( i n x), \ \  n x = n_1 x_1 + n_2 x_2,
$$
are an orthogonal basis in the Hilbert space $H^m$ and, consequently, every function $f \in H^m$ is presented as Fourier series
$$
f = \sum_{n \in {\Bbb Z} \oplus {\Bbb Z}} f_n e^{i n x}
$$
converging to $f$ in the topology of $H^m$.
Thus, one can consider $H^m$, $m \in {\Bbb R}$, as the space of formal Fourier series with the finite norm
$$
\Vert f \Vert^2_m = \sum_n (1 + n \cdot n)^m \cdot |f_n|^2
$$
(see Section 1).
These are the famous Sobolev spaces.

Let $P (x_1, x_2)$ be a homogeneous polynomial of the degree $p \ge 2$ with constant coefficients.
Consider the differential operator $L : H^m \to H^{m - p}$ acting according to the rule
$$
L u = P \bigg( - i \frac{\partial}{\partial x_1}, - i \frac {\partial}{\partial x_2} \bigg) u.
$$

According to [1], define the spaces of infinitely differentiable and generalized functions:
$$
H^\infty = \underset {m}{\cap} H^m, \ \ \,
H^{-\infty} = \underset {m}{\cup} H^m.
$$

\begin{de}
We shall call the operator $L$ {hypoelliptic}, if for any $u \in H^{- \infty}$
the inclusion $L u \in H^\infty$ implies $u \in H^\infty$.
\end{de}


\begin{lem}
\label{Lem1}
The operator $L$ is hypoelliptic if and only if the exist constants $C > 0$ and $k_1$ such that
\begin{equation}
|P(n)| > C (n^2)^{k_1} \ \ \forall \, n.
\label{Pn}
\end{equation}
\end{lem}

{\bf Proof}.
Assume that there exist constants $C > 0$ and $k_1$ such that (5) holds true and $f=\sum_n f_ne^{inx} \in H^\infty$.
Then ${f_n\over P(n)}$ decreases faster than any power of $n$. If there exists a sequence of pairs $n^j$ such that $|P(n^j )| \to 0$ for $j\to\infty$ faster than any power of $|n|$, then, for example, for the functions
$$
f = \sum_j P(n^j )e^{ in^jx} \in H^\infty,
$$
but solution $\sum_j e^{in^jx} \in H^{-2}$, and there is no hypoellipticity.

Recall that in the case of two variables one of necessary and sufficient conditions of hypoellipticity
is the following (see, e.g., [12]):
There exist constants $C$ and $c$ such that
$$
\Bigl| \frac {P^{(\alpha)} (\xi)}{P(\xi)} \Bigr| \le C \, |\xi |^{- |\alpha| c}
$$
for every multi-index $\alpha$ and any $\xi \in {\Bbb R}^2$ large enough.

For $\alpha = p$, the latter inequality gives the inequality $|P(\xi)| \ge C |\xi |^{p c}$, which coincides with inequality \eqref{Pn} on the integer lattice.
This means that every operator hypoelliptic on the plane is hypoelliptic on the torus, but not vice versa.

\begin{pr}
\label{St6}
The operator $L$ is hypoelliptic if and only if for every real root $\alpha$ of the polynomial $P(x, 1)$
there exist constants $C > 0$ and $k$ such that
\begin{equation}
\Bigl| \alpha - \frac pq \Bigr| > \frac C{q^k}
\label{ineq}
\end{equation}
for every rational $p/q$ sufficiently close to $\alpha$.
\end{pr}

{\bf Proof}.
Let $u \in H^{- \infty}$, then $u \in H^m$ with some $m$. Therefore,
$$
u = \sum_n u_n e^{i n x}, \ \   L u = \sum_n u_n P (n)e^{i n x}.
$$

Let $f = \sum_n f_n e^{i n x}$. It is clear that $f \in H^\infty$ if and only if $|f_n| \to 0$ for $n^2 \to \infty$ faster than any power of $n^2$.
Note that if the operator $L$ is hypoelliptic, then the equation $P(n) = 0$ has a unique integer solution $n = 0$.
Indeed, assume that $\nu = (\nu_1, \nu_2)$ is another solution, then the function
$$
{\sum^{+\infty}_{k = -\infty}} e^{i k \nu x} \in H^{- \infty}
$$
belongs to the kernel of $L$, which contradicts the hypoellipticity.

Therefore, the solution of the equation $L u = f$ can be formally written as
$$
u = \sum_{n \not= 0} \frac {f_n}{ P (n)} e^{i n x} + C_0.
$$
The condition $f_0 = 0$ is obviously a condition for the solvability of the equation $L u = f$.
Let us apply Lemma \ref{Lem1}.
Let $\alpha$ be a real root of the polynomial $P(x, 1)$ of multiplicity $r$.
The inequality $|P(n)| > C n^{2k_1}$ is equivalent to
$$
\big | P\big (\tfrac {n_1}{n_2}, 1 \big ) \big | > C (n^2)^{(k_1 - \frac p2)},
$$
since $P \big (\tfrac{n_1}{n_2}, 1 \big)$ tends to zero as $\frac {n_1}{n_2}$ tends to one of the roots.
Thus, we get the inequality
$$
\big | P\big (\tfrac {n_1}{n_2}, 1 \big ) \big | =\big | \tfrac {n_1}{n_2} \alpha \big |^r \big | P_{x_1}^{(r)} \big ( \tfrac {n_1}{n_2}
+ \tau \big ( \alpha - \tfrac{n_1}{n_2} \big ) \big | > C (n^2)^{(k_1 - \tfrac p2)}.
$$

For $\tfrac {n_1}{n_2}$ close enough to $\alpha$, we have
$$
\big | \tfrac {n_1}{n_2} - \alpha \big | > Cn^{2 k_1} > C n_2^{2 k}.
$$
A direct calculation gives $k = \frac 1r \big ( k_1 - \tfrac p2 \big )$.

On the contrary, if $\big | \frac {n_1}{n_2} - \alpha \big | > C n_2^{2 k}$, then $\big | \frac {n_1}{n_2} - \alpha \big | > C n^{2 k} (k < 0)$,
whence we obtain $|P(n)| > C n^{2 k_1}$. The proof is complete.

\medskip

\begin{rem}
{\rm
Inequality \eqref{ineq} is not valid for some transcendental numbers $\alpha$, for example,
$\sum^\infty_{\nu = 1} \frac 1{10^\nu !}$ (see [13]).
On the other hand, the Liouville theorem states that for every algebraic number $\alpha$ of degree $\nu$ inequality \eqref{ineq}
holds true with $k = \nu$ (see [13]).
}
\end{rem}

In particular, we obtain the following

\begin{pr}
\label{St7}
If a polynomial $P$ with integer (or rational) coefficients is irreducible in the field ${\Bbb Q}$, then the operator $L$ is hypoelliptic.
\end{pr}

\medskip

It is easy to see that the hypoelliptic operator $L : H^\infty \to H^\infty$ is reversible.
Here and below, all spaces are assumed to be quotient by the subspace of constants.
The inverse operator $L^{-1}$ acts from $H^m$ into $H^k$ with some $k$.
According to the Thue--Siegel--Roth theorem [13], for every algebraic number $\alpha$ of degree $r \ge 2$ and any
$\varepsilon > 0$ there exists $C > 0$ such that for every rational number $p/q$ the inequality
$$
\Bigl|\alpha - \frac pq \Bigr| > \frac {c}{q^{2 + \varepsilon}}.
$$
holds true.
Then, using calculation of the exponents from the proof of Proposition~\ref{St7}, we obtain the following

\begin{pr}
\label{St8}
Let $r$ be the greatest multiplicity of real roots of the irreducible polynomial $P$.
Then for every $\varepsilon>0$ the operator $L^{-1}$ acts from $H^m$ into $H^{p/2+m-r- \varepsilon}$ continuously.
\end{pr}


\begin{rem}
{\rm
For $p = 2$, by the Liouville theorem, one can put $\varepsilon = 0$.
}
\end{rem}

\end{document}